# One functional property of the $\varsigma$ – function of Riemann

*Azimbay Sadullaev*

National University of Uzbekistan, *sadullaev@mail.ru*

*Abstract.* We prove that if a function $\theta(z) = \int_1^\infty \frac{\pi(t) - Li(t)}{t^{z+1}} dt$, which is holomorphic in $\{\operatorname{Re} z > 1\}$ holomorphically extends to some simply connected domain $G \subset \left\{\operatorname{Re} z > \frac{1}{2}\right\}$, then the $\varsigma(z)$ – function of Riemann has no zeros in this domain, $\varsigma(z) \neq 0 \ \forall z \in G$. As a consequence, it turns out that if the function $\theta(z)$ is holomorphic in $\operatorname{Re} z > \frac{1}{2}$, then the Riemann hypothesis has a positive solution.

**0. Asknowledjment.** I have never dealt with the Riemann $\varsigma$ – function, although I heard about the difficult Riemann hypothesis about this function. In 2020, my friend Alexander Kytmanov, an excellent mathematician, approached me and asked me to help in some applications of entire functions to functional properties of $\varsigma$ – function. As a result, we published a joint article: *A. Kytmanov, V. Kuzovatov, A. Sadullaev,* On the Zeta-Function of Zeros of Entire Function, J. Sib. Fed. Univ. Math. Phys., 2021, 14(5), 599–603. In the process of studying the functional properties of $\varsigma$ – function, I also independently proved a number of statements, related to the Riemann hypothesis about the zeros of the Riemann function. In this article, we will talk about this function and prove some functional properties of the $\varsigma$ – function and its zeros.

**1. Introduction.** Function

$$\varsigma(x) = \sum_{k=1}^\infty \frac{1}{k^x}, \ x \in \mathbb{R}, \qquad (1)$$

is called the Riemann zeta-function. It is clear that for $x > 1$ the series converges, and for $x \leq 1$ it diverges. We define the zeta function for complex numbers $z = x + iy \in \mathbb{C}$ by setting

$$\varsigma(z) = \sum_{k=1}^\infty \frac{1}{k^z} = \sum_{k=1}^\infty \frac{1}{e^{\ln k^z}} = \sum_{k=1}^\infty \frac{1}{e^{z \ln k}}.$$

We see that the series converges uniformly on any compact set of the domain $\{\operatorname{Re} z > 1\}$ and, therefore, its sum is analytic in this domain $\{\operatorname{Re} z > 1\}$. This function plays a big role in number theory. The great mathematician Euler spent a lot of time studying this function and proved a number of identities and properties. In particular:

a). In the domain $\{\operatorname{Re} z > 1\}$, $\varsigma$ – function is represented as

$$\varsigma(z) = \prod_p \left( \frac{1}{1 - p^{-z}} \right), \quad p - \text{prime numbers.} \tag{2}$$

This implies that $\varsigma$ – function has no zeros in $\{\operatorname{Re} z > 1\}$.

b). The Riemann $\varsigma(z)$ – function extends holomorphically to the entire plane, except the point $z = 1$, where it has a simple pole: at $\{\operatorname{Re} z > 0\}$ it is represented as

$$\varsigma(z) = \frac{1}{z-1} + \frac{1}{2} + z \int_1^\infty \frac{\rho(t)}{t^{z+1}} dt, \quad \rho(t) = 1 - \{t\}. \tag{3}$$

c). The zeros of the Riemann function on $\{\operatorname{Re} z \leq 0\}$ lie only at points $-2, -4, \ldots$, (trivial zeros) and on $\{\operatorname{Re} z \geq 1\}$ it has no zeros. We are required to study them in the the band $\{0 < \operatorname{Re} z < 1\}$: **prove or disprove that all non-trivial zeros of $\varsigma$ – function (in $\{0 < \operatorname{Re} z < 1\}$) lie on the vertical line $\left\{\operatorname{Re} z = \frac{1}{2}\right\}$!**

d). The non-trivial zeros of the function are symmetric with respect to the axis and the axis $\left\{\operatorname{Re} z = \frac{1}{2}\right\}$. They do not lie on vertical lines $\{\operatorname{Re} z = 0\}$, $\{\operatorname{Re} z = 1\}$ and in the interval $\{0 < \operatorname{Re} z < 1, \operatorname{Im} z = 0\}$. There are infinitely many of them, more than 10 zeros are calculated on the computer, they all look like $z = \frac{1}{2} + iy$. On the upper plane, the first zeros are equal:

$1\backslash 2 + i \cdot 14.134\ldots$; $1\backslash 2 + i \cdot 21.022\ldots$; $1\backslash 2 + i \cdot 25.010\ldots$; $1\backslash 2 + i \cdot 30.424\ldots$;
$1\backslash 2 + i32.935\ldots$; $1\backslash 2 + i37.935\ldots$; $1\backslash 2 + i \cdot 40.918\ldots$; $1\backslash 2 + i \cdot 43.327\ldots$;
$1\backslash 2 + i \cdot 48.005$; $1\backslash 2 + i \cdot 49.773\ldots$.

According to (3), $(z-1)\varsigma(z) = 1 + (z-1)\left[\frac{1}{2} + z\int_1^\infty \frac{\rho(t)}{t^{z+1}} dt\right] = 1 + \psi(z)$ and

$$\frac{\varsigma'(z)}{\varsigma(z)} = -\frac{1}{z-1} + \frac{\psi'(z)}{(z-1)\varsigma(z)} = -\frac{1}{z-1} + F(z), \quad F(z) = \frac{\psi'(z)}{(z-1)\varsigma(z)} \in O(\Pi). \tag{4}$$

Here $\Pi = \{\operatorname{Re} z > 1/2, \} \cap \{|\operatorname{Im} z| < 4\pi\}$.

**2. Analytic properties of $\varsigma$ – function.** First, let $\{\operatorname{Re} z > 1\}$. Then

$$\varsigma(z) = \exp\left\{-\ln\prod_p\left(1-p^{-z}\right)\right\} = \exp\left\{\sum_p p^{-z} + \sum_p\left[\frac{p^{-2z}}{2} + \frac{p^{-3z}}{3} + \ldots\right]\right\} =$$

$$= \exp\left\{\sum_p p^{-z} + f(z)\right\}, \quad f(z) = \sum_p\left[\frac{p^{-2z}}{2} + \frac{p^{-3z}}{3} + \ldots\right] \in O\{\operatorname{Re} z > 1/2\};$$

$$\varsigma(z) = \exp\left[\sum_p p^{-z}\right]\exp[f(z)].$$

$$\exp\left[\sum_p p^{-z}\right] = \varsigma(z)\exp[-f(z)] \in O\{\operatorname{Re} z > 1/2\} \setminus \{z=1\};$$

(5)

$$(z-1)\exp\left[\sum_p p^{-z}\right] = (z-1)\varsigma(z)\exp[-f(z)] \in O\{\operatorname{Re} z > 1/2\}$$

and $\neq 0$ in $\Pi = \{\operatorname{Re} z > 1/2,\}\cap\{|\operatorname{Im} z| < 4\pi\}$.

Note that the series $\sum_p p^{-z}$ converges uniformly inside the domain $\{\operatorname{Re} z > 1\}$

and $\dfrac{d}{dz}\sum_p p^{-z} = -\sum_p p^{-z}\ln p$.

**3. Expression for** $\dfrac{d}{dz}\sum_p p^{-z} = -\sum_p p^{-z}\ln p$. We take

$\varphi(t) = \pi(t) - Li(t), \, d\varphi(t) = \sum_p \delta_p - \dfrac{1}{\ln t}$, where $\delta_p$ is the Dirac function at the point $p$, $\pi(t)$ is the number of prime numbers in the segment $[2,t]$, and $Li(t) = \displaystyle\int_2^t \dfrac{ds}{\ln s}$.

For $\operatorname{Re} z > 1$

$$\int_2^\infty t^{-z}\ln t\, d\varphi(t) = t^{-z}\ln t\,\varphi(t)\Big|_2^\infty - \int_2^\infty\left[-zt^{-z-1}\ln t + t^{-z-1}\right]\varphi(t)dt =$$

$$= -\frac{\varphi(2)\ln 2}{2^z} - \int_2^\infty\left[-zt^{-z-1}\ln t + t^{-z-1}\right]\varphi(t)dt.$$

$$\left\{u = t^{-z}\ln t,\; du = \left[-zt^{-z-1}\ln t + t^{-z-1}\right]dt,\; dv = d\varphi(t),\; v = \varphi(t)\right\}$$

On the other side,

$$\int_2^\infty t^{-z}\ln t\, d\varphi(t) = \sum_p p^{-z}\ln p - \int_2^\infty t^{-z}\ln t \frac{1}{\ln t}dt = \sum_p p^{-z}\ln p - \frac{2^{1-z}}{z-1}.$$

From here,

$$\sum_p p^{-z}\ln p = \frac{2^{1-z}}{z-1} - \frac{\varphi(2)\ln 2}{2^z} + \int_2^\infty \left[zt^{-z-1}\ln t - t^{-z-1}\right]\varphi(t)dt. \qquad (6)$$

**4. Function** $\Phi(z) = \int_2^\infty \left[zt^{-z-1}\ln t - t^{-z-1}\right]\varphi(t)dt$. Differentiating

$$\exp\left[\sum_p p^{-z}\right] = \varsigma(z)\exp[-f(z)]$$

in $\{\mathrm{Re}\, z > 1\}$ we get

$$-\exp\left[\sum_p p^{-z}\right]\sum_p p^{-z}\ln p = \varsigma'(z)\exp[-f(z)] - \varsigma(z)\exp[-f(z)]f'(z);$$

$$-\varsigma(z)\sum_p p^{-z}\ln p = \varsigma'(z) - \varsigma(z)f'(z);$$

$$\sum_p p^{-z}\ln p = f'(z) - \frac{\varsigma'(z)}{\varsigma(z)}$$

and according to (4)

$$\sum_p p^{-z}\ln p = f'(z) - \frac{\varsigma'(z)}{\varsigma(z)} = \frac{1}{z-1} + f'(z) - F(z); \quad F(z) \in O(\Pi).$$

From here,

$$\Phi(z) = \int_2^\infty \left[zt^{-z-1}\ln t - t^{-z-1}\right]\varphi(t)dt = \frac{1}{z-1} - \frac{2^{1-z}}{z-1} + \frac{\varphi(2)\ln 2}{2^z} + f'(z) - F(z) =$$

$$= \frac{1}{z-1} - \frac{e^{-(z-1)\ln 2}}{z-1} + \frac{\varphi(2)\ln 2}{2^z} + f'(z) - F(z) \in O(\Pi), \qquad (7)$$

since $\dfrac{1}{z-1} - \dfrac{e^{-(z-1)\ln 2}}{z-1} \in O(\mathbb{C})$.

**5. Expansion in a series.** We have

$$\Phi(z) = \int_2^\infty \left[-zt^{-z-1}\ln t + t^{-z-1}\right]\varphi(t)dt = \left\{\int_2^\infty \frac{n(t) - Li(t)}{t^{z+1}}dt - z\int_2^\infty \ln t \frac{n(t) - Li(t)}{t^{z+1}}dt\right\} \in O(\Pi);$$

$$\theta(z) = \int_2^\infty \frac{n(t) - Li(t)}{t^{z+1}}dt, \quad \theta^{(n)}(a) = (-1)^n \int_2^\infty \frac{n(t) - Li(t)}{t^{a+1}}\ln^n t\, dt, \quad \mathrm{Re}\, a > 1.$$

Then $\Phi(z) = \theta(z) + z\theta'(z) = \theta(z) + (z-1)\theta'(z) + \theta'(z)$ and it is holomorphic in $\{\operatorname{Re} z > 1\}$, because, $\pi(t) - Li(t) = O\left(t \cdot \exp\left[-0,005(\ln t)^{\frac{3}{5}-0}\right]\right)$, therefore, asymptotically $\dfrac{\pi(t) - Li(t)}{t} > \dfrac{1}{\ln^\delta t} \ \forall \delta > 0$.

We fix $0 < \varepsilon < 1/2$, $|b| < 3\pi$, $a = 1 + \varepsilon + ib$. Expanding into a series, we get

$$\theta(z) = \theta(a) + \frac{\theta'(a)}{1!}(z-a) + \ldots + \frac{\theta^{(n)}(a)}{n!}(z-a)^n + \ldots$$

$$\Phi(z) = \theta(z) + (z-a)\theta'(z) + \theta'(z) = \sum_{n=0}^{\infty} \frac{\theta^{(n)}(a)}{n!}(z-a)^n + \sum_{n=1}^{\infty} \frac{\theta^{(n)}(a)}{(n-1)!}(z-a)^n +$$

$$+ \sum_{n=1}^{\infty} \frac{\theta^{(n)}(a)}{(n-1)!}(z-a)^{n-1} = \theta(a) + \sum_{n=1}^{\infty}\left[\frac{\theta^{(n)}(a)}{n!} + \frac{\theta^{(n)}(a)}{(n-1)!} + \frac{\theta^{(n+1)}(a)}{n!}\right](z-a)^n =$$

$$= \theta(a) + \sum_{n=1}^{\infty} \frac{n+1}{n!}\left[\theta^{(n)}(a) + \frac{\theta^{(n+1)}(a)}{(n+1)}\right](z-a)^n;$$

Since $\Phi(z) \in O(\Pi)$, then $\overline{\lim_{n\to\infty}}\left|\dfrac{1}{n!}\left[\theta^{(n)}(a) + \dfrac{\theta^{(n+1)}(a)}{(n+1)}\right]\right|^{1/n} \leq \dfrac{1}{2}$.

**6. The reasoning on this point was prompted to me by my student, Professor Isroil Ikramov. I am very grateful to him.** Let

$$g(z) = \sum_{n=1}^{\infty} \frac{\theta^{(n)}(a)}{n!}(z-a)^n \text{ and } \overline{\lim_{n\to\infty}}\left|\frac{1}{n!}\left[\theta^{(n)}(a) + \frac{\theta^{(n+1)}(a)}{(n+1)}\right]\right|^{1/n} \leq \frac{1}{2}.$$

Then

$$\sum_{n=1}^{\infty} \frac{1}{n!}\left[\theta^{(n)}(a) + \frac{\theta^{(n+1)}(a)}{(n+1)}\right](z-a)^n = g(z) + \frac{g(z)}{z-a} - \theta(a) = H(z) \in O\left\{|z-a| < \frac{1}{2}\right\}.$$

$$g(z) = \frac{(z-a)[H(z) + \theta(a)]}{z - (\varepsilon + ib)} \in O\left\{|z-a| < \frac{1}{2}\right\}.$$

Hence, $\overline{\lim_{n\to\infty}}\left|\dfrac{\theta^{(n)}(a)}{n!}\right|^{1/n} \leq \dfrac{1}{2}$ and the series

$\theta(z) = \theta(a) + \dfrac{\theta'(a)}{1!}(z-a) + \ldots + \dfrac{\theta^{(n)}(a)}{n!}(z-a)^n + \ldots$ converges uniformly at

$|z-a|<1/2$, its sum is holomorphic there. Hence, a holomorphic function $\theta(z) = \int_{2}^{\infty} \frac{\pi(t)-Li(t)}{t^{z+1}} dt$ in $\{\operatorname{Re} z > 1\}$ holomorphically extends to $\{\operatorname{Re} z > 1\} \cup \hat{\Pi}$, $\hat{\Pi} = \{\operatorname{Re} z > \frac{1}{2},\ |\operatorname{Im} z| < 3\pi\}$.

**7. Theorem.** If the function $\theta(z)$ extends holomorphically to some simply connected domain $G \subset \{\operatorname{Re} z > \frac{1}{2}\}$, then the Riemann $\varsigma(z)$-function has no zeros in this domain.

**Corollary.** If the holomorphic in $\{\operatorname{Re} z > 1\}$ function $\theta(z) = \int_{2}^{\infty} \frac{\pi(t)-Li(t)}{t^{z+1}} dt$ holomorphically extends to $\{\operatorname{Re} z > \frac{1}{2}\}$, then all nontrivial zeros of the Riemann function lie on the vertical line $\{\operatorname{Re} z = \frac{1}{2}\}$.

P r o o f   of   T h e o r e m.   According to (6)
$$\sum_{p} p^{-z} \ln p = \frac{2^{1-z}}{z-1} - \frac{\varphi(2)\ln 2}{2^z} + \int_{2}^{\infty} \left[ zt^{-z-1}\ln t - t^{-z-1} \right] \varphi(t) dt = \frac{2^{1-z}}{z-1} - \frac{\varphi(2)\ln 2}{2^z} + \Phi(z).$$

Therefore, the function $\sum_{p} p^{-z} \ln p$ extends holomorphically to a simply connected domain $G \setminus \left[\frac{1}{2}, 1\right]$ and $\sum_{p} p^{-z} = \int_{[2,z]} \sum_{p} p^{-\xi} \ln p - \sum_{p} p^{-2}$, where $[2,z]$ is a smooth curve connecting points $2, z$ inside a simply connected domain $G \setminus \left[\frac{1}{2}, 1\right]$. It follows, that $\sum_{p} p^{-z} \ln p$ also extends holomorphically to $G \setminus \left[\frac{1}{2}, 1\right]$. From here and from $\varsigma(z) = \exp\left[ \sum_{p} p^{-z} \right] \exp[f(z)]$ it follows that $\varsigma(z) \neq 0$ in the domain $G \setminus \left[\frac{1}{2}, 1\right]$.

Since the function $\varsigma(z) \neq 0$ on the interval $\{0 < \operatorname{Re} z < 1,\ \operatorname{Im} z = 0\}$, then it is $\neq 0$ in $G$. *The theorem is proven.*

**Task 1.** Using the expansion of a function $\theta(z) = \int_2^\infty \frac{\pi(t) - Li(t)}{t^{z+1}} dt$ in a Taylor series, $\theta(z) = \theta(a) + \frac{\theta'(a)}{1!}(z-a) + \ldots + \frac{\theta^{(n)}(a)}{n!}(z-a)^n + \ldots$ we have proved that the radius of convergence of this series $R(a) \geq \frac{1}{2}$ for all $a = 1 + \varepsilon + ib$, $|b| < 3\pi$. Prove that $R(a) \geq \frac{1}{2} \ \forall b \in \mathbb{R}_+$. The positive answer to this task is equivalent to the Riemann hypothesis.

**More general task 2.** Let a function $\sigma(t)$, $t \in [1, \infty)$, have the function property of $\frac{\pi(t) - Li(t)}{t}$; that is sign-variable and asymptotically $|\sigma(t)| < \ln^\delta t, \ \forall \delta > 0$. It follows that the integral $\phi(z) = \int_1^\infty \frac{\sigma(t)}{t^z} dt$ converges absolutely and uniformly in a closed domain $\operatorname{Re} z \geq 1$. Let moreover, it is known that the function $\phi(z)$ extends holomorphically to $\hat{\Pi} = \left\{ \operatorname{Re} z > \frac{1}{2}, \ |\operatorname{Im} z| < 3\pi \right\}$. Prove that it extends holomorphically to $\operatorname{Re} z > \frac{1}{2}$.

**Example** (I. Ikromov, counterexample to task 2). Consider the function $\sigma(t) = \frac{2\cos(12 \ln t)}{t^\gamma}$, $t \in [1, \infty)$, $\gamma < \frac{1}{4}$.

By its asymptotic nature, this function is close to $\frac{\pi(t) - Li(t)}{t}$, in the sense that for $t \to \infty$, they are sign-variable and tend to zero faster than $\frac{1}{\ln^\delta t} \ \forall \delta > 0$. Let

$$\phi(z) = \int_1^\infty \frac{\sigma(t)}{t^z} dt. \tag{9}$$

The integral converges at $\operatorname{Re} z > 1 - \delta$ and is holomorphic in the domain $\operatorname{Re} z > 1 - \delta$. It can be explicitly calculated by the change $e^\tau = t$

$$\phi(z) = \int_1^\infty \frac{\sigma(t)}{t^z} dt = \int_0^\infty \frac{2\cos 12\tau}{e^{\tau(z+\gamma)}} e^\tau d\tau = \frac{2(z+\gamma-1)}{12^2 + (1-z-\gamma)^2}.$$

Thus, $\phi(z) = \int_1^\infty \frac{\sigma(t)}{t^z} dt$ is holomorphic in $\mathbb{C} \setminus \{1 - \gamma \pm 12i\}$.

This example shows that holomorphy of $\theta(z) = \int_2^\infty \frac{\pi(t) - Li(t)}{t^{z+1}} dt$ in $\left\{ \operatorname{Re} z > \frac{1}{2} \right\}$ does not follow from the elementary properties of the function $\frac{\pi(t) - Li(t)}{t}$. Since the holomorphy of $\theta(z) \in O\left( \left\{ \operatorname{Re} z > \frac{1}{2} \right\} \right)$ is equivalent to the Riemann hypothesis, the proof of the Riemann hypotheses requires additional finer properties of the function $\frac{\pi(t) - Li(t)}{t}$.